\documentclass[a4paper]{amsart}
\usepackage[latin1]{inputenc}
\usepackage[all]{xy}
\usepackage{amssymb}
\usepackage{rotating}
\usepackage{amsmath}
\usepackage{graphicx}
\usepackage{epsfig}
\usepackage{eufrak}

\newtheorem{theorem}{Theorem}[section]
\newtheorem{rem}[theorem]{Remark}
\newtheorem{prop}[theorem]{Proposition}

\newtheorem{cor}[theorem]{Corollary}
\newtheorem{defi}{Definition}[section]

\newcommand{\bprf}{{\it Proof.~}}

\newcommand{\ra}{\rightarrow}
\newcommand{\eprf}{\hfill $\square$ \smallskip\par}
\newcommand{\erem}{\hfill $\square$}

\newcommand{\C }{ \mathbb{C}}

\newcommand{\Z}{\mathbb{Z}}
\newcommand{\Q}{\mathbb{Q}}

\newcommand{\Cf}{C_{f_g}}
\newcommand{\Cfu}{C_{f_1}}
\newcommand{\Cfd}{C_{f_2}}
\newcommand{\Cft}{C_{f_3}}
\newcommand{\Yf}{Y_{f_g}}
\newcommand{\Yfd}{Y_{f_2}}
\newcommand{\Zf}{Z_{f_g}}
\newcommand{\Sf}{S_{f_g}}
\newcommand{\Vf}{V_{f_g}}
\newcommand{\Sfd}{S_{f_2}}
\newcommand{\Sft}{S_{f_3}}
\newcommand{\pc}{P_i\times C_j}
\newcommand{\pd}{P_i\times D_j}
\newcommand{\qc}{Q_i\times C_j}
\newcommand{\qd}{Q_i\times D_j}

\newcommand{\rmd}{\mbox{d}}
\newcommand{\norm}[1]{\left|\left|#1\right|\right|}

\def\DynkinDD#1#2#3#4#5
   {\begin{picture}(124, 70)%
   \put(5, 40){\circle{4}}%
   \put(7, 40){\line(1, 0){28}}%
   \put(37,40){\circle{4}}%
   \put(39, 40){\line(1, 0){28}}%
   \put(69,40){\circle{4}}%
   \put(37,38){\line(0,-11){28}}%
   \put(37,42){\line(0,11){28}}%
   \put(37,72){\circle{4}}%
   \put(37,8){\circle{4}}%
   \put(5, 46){\makebox(0, 0)[b]{$#1$}}%
   \put(45,8){\makebox(0, 0)[b]{$#2$}}%
   \put(45,46){\makebox(0, 0)[b]{$#3$}}%
   \put(45,62){\makebox(0, 0)[b]{$#4$}}%
   \put(69, 16){\makebox(0, 0)[b]{$#5$}}%
  \end{picture}}

\def\DynkinEIII#1#2#3#4#5#6#7#8
   {\begin{picture}(240, 65)%
   \put(5, 45){\circle{4}}%
   \put(7, 45){\line(1, 0){28}}%
   \put(37, 45){\circle{4}}%
   \put(39, 45){\line(1, 0){28}}%
   \put(69, 45){\circle{4}}%
   \put(71, 45){\line(1, 0){28}}%
   \put(101, 45){\circle{4}}%
   \put(103, 45){\line(1, 0){28}}%
   \put(133, 45){\circle{4}}%
   \put(135, 45){\line(1, 0){28}}%
   \put(165, 45){\circle{4}}%
   \put(167, 45){\line(1, 0){28}}%
   \put(197, 45){\circle{4}}%
   \put(101, 43){\line(0, -11){28}}%
   \put(101, 13){\circle{4}}%
   \put(5, 51){\makebox(0, 0)[b]{$#1$}}%
   \put(37, 51){\makebox(0, 0)[b]{$#2$}}%
   \put(69, 51){\makebox(0, 0)[b]{$#3$}}%
   \put(101, 51){\makebox(0, 0)[b]{$#4$}}%
   \put(133, 51){\makebox(0, 0)[b]{$#5$}}%
   \put(165, 51){\makebox(0, 0)[b]{$#6$}}%
   \put(193, 51){\makebox(0, 0)[b]{$#7$}}%
    \put(105, 13){\makebox(0, 0)[lb]{$#8$}}%
\end{picture}}

\input epsf
\input cyracc.def

\makeatletter
\def\blfootnote{\xdef\@thefnmark{}\@footnotetext}

\begin{document}

\title{New families of Calabi--Yau 3--folds without maximal unipotent monodromy}
\author{Alice Garbagnati}
\address{Alice Garbagnati, Dipartimento di Matematica, Universit\`a di Milano,
  via Saldini 50, I-20133 Milano, Italia}

\email{alice.garbagnati@unimi.it}

\begin{abstract}
The aim of this paper is to construct families of Calabi--Yau 3-folds without boundary points with maximal unipotent monodromy and to describe the variation of their Hodge structures. In particular five families are constructed. In all these cases the variation of the Hodge structures of the Calabi--Yau 3-folds is basically the variation of the Hodge structures of a family of curves. This allows us to write explicitly the Picard--Fuchs equation for the 1-dimensional families.\\
These Calabi--Yau 3-folds are desingularizations of quotients of the product of a (fixed) elliptic curve and a K3 surface admitting an automorphisms of order 4 (with some particular properties). We show that these K3 surfaces admit an isotrivial elliptic fibration.    
\end{abstract}

\maketitle

\markboth{FAMILIES OF CALABI--YAU 3--FOLDS WITHOUT MAXIMAL UNIPOTENT MONODROMY}{ALICE GARBAGNATI}

\blfootnote {{\it 2000 Mathematics Subject Classification: 14J32, 14J28, 14J50}.} \blfootnote {{\it Key words: Calabi--Yau, automorphisms, K3 surfaces, Picard--Fuchs equations}.}
\section{Introduction}

A Calabi--Yau variety $Y$ is a smooth compact K\"ahler manifold with a trivial canonical bundle and such that $H^i(Y,\mathcal{O}_Y)=0$ if $0<i<dim (Y)$. The main conjecture on the families of Calabi--Yau 3--folds is the mirror symmetry conjecture that, roughly speaking, predicts the existence of a mirror family for each family of Calabi--Yau 3-folds. This conjecture is stated for families of Calabi--Yau 3--folds admitting a boundary point with maximal unipotent monodromy, otherwise it is not clear how to define the mirror family. For this reason it seems important to construct examples of families of Calabi--Yau 3-folds without maximal unipotent monodromy, to describe the variation of their Hodge structures and their parametrizing space and to give explicitly their Picard--Fuchs equations.\\

We recall that an automorphism on a Calabi--Yau $Y$ of dimension $n$ is symplectic if it acts as the identity on the 1-dimensional space $H^{n,0}(Y)$, otherwise it is non symplectic. In order to construct a Calabi--Yau of a given dimension one can consider the desingularization of certain quotients of products of Calabi--Yau varieties of lower dimension by finite groups of automorphisms. This construction is very well known (\cite{borcea}, \cite{cynk hulek}, \cite{Rohde thesis}, \cite{Rohde}, \cite{voisin miroirs}). In particular Rohde (\cite{Rohde}) produces in this way families of Calabi--Yau 3-folds parametrized by Shimura varieties. The Calabi--Yau 3-folds considered in \cite{Rohde} are obtained as quotient of the product of a K3 surface and an elliptic curve by the group $\Z/3\Z$, and both the K3 surface and the elliptic curve considered admit $\Z/3\Z$ as group of non symplectic automorphisms. In \cite{noi} the K3 surfaces used in the Rohde's construction are shown to be desingularizations of the quotient of the product of two curves by the group $\Z/3\Z$ and hence the Calabi--Yau 3-folds are desingularizations of the quotient of the product of three curves by the group $(\Z/3\Z)^2$. This allowed us to determine explicitly the Picard--Fuchs equation of a 1-dimensional family.\\

Here we use a similar construction, considering K3 surfaces and elliptic curves admitting non symplectic automorphisms of order 4. The main problem is that there is not yet a complete classification of the K3 surfaces admitting a non symplectic automorphism of order 4. However the construction used in \cite{noi} can be considered also in this situation. Indeed we will consider the product of the elliptic curve $E_i$ (i.e.\ the elliptic curve admitting a non symplectic automorphism $\alpha_E$ of order 4 ) with a curve $\Cf$ of genus $g$ which admits an automorphism $\alpha_C$ acting with eigenvalues $i$ and $-i$ on the space of the holomorphic 1-forms and such that one of the eigenspaces has dimension 1. The quotient of the product $E_i\times \Cf$ by an automorphism obtained from $\alpha_E$ and $\alpha_C$ admits a desingularization which is a K3 surface $\Sf$ with a non symplectic automorphism, $\alpha_S$, of order 4.\\

In sections \ref{sec: constructing K3 Sd},  \ref{subsection: isotrivial}, \ref{section: geometrical description} we construct K3 surfaces with a non symplectic automorphism of order 4 such that its fixed locus and the fixed locus of its square consist of isolated points and rational curves. In particular in Section \ref{sec: constructing K3 Sd} we construct the K3 surfaces $\Sf$ as quotient of the product of an elliptic curve $E_i$ and a curve $\Cf$ (described in Section \ref{sec: the curves Cf}) by the group $\Z/4\Z$. In Section \ref{subsection: isotrivial} we describe the family $Q_{b^2}$ of K3 surfaces as isotrivial fibration with section. In Section \ref{section: geometrical description} we observe that the surfaces $\Sf$ admit a description as $4:1$ cover of $\mathbb{P}^2$ and we construct another family, $R$, of such a K3 surfaces.\\
In Section \ref{section: Calabi Yau} we consider the quotient of $E_i\times S$ by the automorphism $\alpha_E^3\times\alpha_S$, where $S$ is one of the K3 surfaces constructed. This quotient admits a desingularization which is a Calabi--Yau 3-fold. The main result of the Section \ref{section: Calabi Yau} is the description of the Hodge structure of this Calabi--Yau 3-fold and of its variation. 
In the Section \ref{section: PF} we show that the Picard--Fuchs equations of the families have no a solution with maximal unipotent monodromy. We write these equations explicitly for the two 1-dimensional families we constructed.\\
The properties of the families of Calabi--Yau 3-folds obtained and of the K3 surfaces and curves involved in their construction are summarized in the following table:
$$
\begin{array}{|c|c|c|c|c|c|c|}\hline
\mbox{Family of}&h^{1,1}&h^{2,1}&\mbox{Related family of}&\mbox{Transcendental}&\mbox{Related family}&\mbox{Parametrizing}\\
\mbox{Calabi--Yau's}&&&\mbox{K3 surfaces}&\mbox{lattice of K3}&\mbox{of curves}&\mbox{spaces}\\
\hline
Y_{SI}&90&0&X_{SI}&U(2)&C_{f_1}\simeq E_i&point\\
\mbox{(Remark }\ref{rem: YSI}\mbox{)}&&&\mbox{(Section }\ref{section: Shioda Inose}\mbox{)}&&\mbox{(Section }\ref{sec: the curves Cf}\mbox{)}&\mbox{(Corollary }\ref{cor: moduli space}\mbox{)}\\
\hline
Y_{f_2}&73&1&S_{f_2}&\langle 2 \rangle^2\oplus \langle 2\rangle\oplus\langle -2\rangle&C_{f_2}&\mathbb{B}^1\\
\mbox{(Proposition }\ref{prop: Hodge numbers CY}\mbox{)}&&&\mbox{(Section }\ref{prop: elliptic fibration on E}\mbox{)}&&\mbox{(Section }\ref{sec: the curves Cf}\mbox{)}&\mbox{(Corollary }\ref{cor: moduli space}\mbox{)}\\
\hline
Y_{f_3}&56&2&S_{f_3}&U(2)^2\oplus\langle -2\rangle^2&C_{f_3}&\mathbb{B}^2\\
\mbox{(Proposition }\ref{prop: Hodge numbers CY}\mbox{)}&&&\mbox{(Section }\ref{prop: elliptic fibration on E}\mbox{)}&&\mbox{(Section }\ref{sec: the curves Cf}\mbox{)}&\mbox{(Corollary }\ref{cor: moduli space}\mbox{)}\\
\hline
W_{b^2}&61&1&Q_{b^2}&U(2)^2&E&\mathbb{H}^1\\
\mbox{(Section }\ref{subsection: Wb}\mbox{)}&&&\mbox{(Proposition }\ref{prop: Qb2}\mbox{)}&&\mbox{(elliptic curve)}&\mbox{(Section }\ref{subsection: Wb}\mbox{)}\\
\hline
M&39&3&R&U(2)^2\oplus \langle-2\rangle^4&&\\
\mbox{(Section }\ref{subsection: 3 dim fam CY}\mbox{)}&&&\mbox{(Proposition }\ref{prop: a 3 dimensional family}\mbox{)}&&&\\
\hline
\end{array}
$$  $ $\\
{\it I'm grateful to Bert van Geemen for useful suggestions and for his kind support during the preparation of this paper.}

\section{The curves $\Cf$}\label{sec: the curves Cf}
In order to construct K3 surfaces and Calabi--Yau 3-folds with a non symplectic automorphism of order 4, we introduce the curves $\Cf$. These curves admit a particular automorphism of order 4 and our construction can be done only for certain curves of genus $g=1,2,3$.\\

Let us consider the hyperelliptic curves of genus $g$
\begin{equation}\label{equation: Cf}\Cf:\ \ z^2=rf_g(r^2),\ \ \ \ \ deg(f_g)=g,\ \ \ \ \ f_g\mbox{ without multiple roots.}\end{equation}
A basis of the space of the 1-holomorphic forms on $\Cf$ is
given by $\{r^a dr/z\}_{0\leq a\leq g-1}$. Set $\eta_C:=dr/z$, $\omega_C:=rdr/z$ and
$\nu_C:=r^2dr/z$.\\
The curve $\Cf$ admits the automorphism $\alpha_C:(r,z)\ra
(-r,iz)$. \\
We have
$\alpha_C^*(\omega_C)=-i\omega_C$, $\alpha_C^*(\eta_C)=i\eta_C$ and
$\alpha_C^*(\nu_C)=i\nu_C$.
In particular the action of $\alpha_C$ on $H^{1,0}(\Cf)$ has only $i$ and $-i$ as eigenvalues and one of them has an eigenspace of dimension 1. By \cite[Lemma 7.6]{bert picard} the curves admitting an automorphism $\alpha_C$ of order 4 such that its action on the spaces of the one holomorphic form has only $i$ and $-i$ as eigenvalues and one eigenspace has dimension 1, are exactly the curves $\Cf$ given in \eqref{equation: Cf}. \\
We observe that if $g=3$, the eigenspace of $H^{1,0}(\Cf)$ of dimension 1 is relative to the eigenvalue $-i$ and if $g=1$,
the eigenspace of $H^{1,0}(\Cf)$ of dimension 1 is relative to the eigenvalue $i$. For this reason in the following these two cases are analyzed separately (cf.\ Sections \ref{section: Shioda Inose}, \ref{construction K3}).\\
One can give a different equation for the curves $\Cf$, which exhibits them as 4:1 cover of $\mathbb{P}^1$. Indeed let us consider the curves \begin{equation}\label{equation D}D_g:\ \ w^4=h_2(t)l_{g}(t)^2,\ \ \ g=1,2,3,\ \ \deg(h_2)=2,\ \ \deg(l_{g})=g,\ \ h_2l_g\mbox{ without multiple roots.}\end{equation}
These curves are 4:1 cover of $\mathbb{P}^1$ branched over the zeroes of $h_2$ and  of $l_{g}$. By Riemann--Hurwitz formula, the genus of $D_g$ is $g$. Of course these curves admit the order 4 automorphism $\gamma_D:(w,t)\mapsto (iw, t)$ which is the cover automorphism.\\
In \cite[Remark 7.7]{bert picard} it is shown that the curves $C_{f_3}$ and $D_3$ are isomorphic. Let us focus on case $g=2$. 
Up to a choice of coordinate of $\mathbb{P}^1$ one can assume that the equation of the curve $D_2$ is $w^4=t(t-1)^2(t-\lambda)^2$. The map $(w,t)\ra r:=w^2/((t-1)(t-\lambda))$ is a $2:1$ map  $D_2\ra \mathbb{P}_r^1$ (this can be easily check by writing the divisor associated to the meromorphic function $w^2/((t-1)(t-\lambda))$). By the definition of $r$, we find $w^2=r(t-1)(t-\lambda)$ and substituting the new variable $r$ in the equation of $D_2$, we have $r^2=t$. Thus we find a new equation for $D_2$, $w^2=r(r^2-1)(r^2-\lambda)$, which is, up to a choice of coordinates of $\mathbb{P}^1$, the equation of $C_{f_2}$. Moreover the cover automorphism $\gamma_D$ induces the automorphism $(w,r)\mapsto (iw,-r)$, which is the automorphism $\alpha_C$ on the curve $C_{f_2}$.\\
The case $g=1$ is very similar. Hence we proved that the curves $\Cf$ endowed with the automorphism $\alpha_C$ coincide with the curves $D_g$ endowed with the automorphism $\gamma_D$.\\

We recall the following result on moduli space of the curves $\Cf$:
\begin{prop}{\rm (cf.\ \cite{bert picard})}\label{prop: moduli space Cf} The curves $\Cf$ are parametrized by the symmetric domain $U(1,g-1)/(U(1)\times U(g-1))\simeq \mathbb{B}^{g-1}:=\{w\in \C^{g-1}| \norm{w}\leq 1\}$ .\end{prop}

\section{The K3 surfaces  $\Sf$ with a non symplectic automorphism of order 4.}\label{sec: constructing K3 Sd}
In this section we construct the K3 surfaces $\Sf$, which are the desingularization of the quotient of product of the elliptic curve $E_i$ (the elliptic curve with a non symplectic automorphism of order 4) and the curve $\Cf$ by an automorphism of order 4. We first consider the case $g=1$ and after that the other two cases because the choice of the automorphism of order 4 depends on the eigenvalues relative to the eigenspace with dimension 1. 
\subsection{Case $g=1$: a K3 surface described by Shioda and Inose}\label{section: Shioda Inose}
The curve $\Cf$, in case $g=1$, is the elliptic curve $E_i$. We recall a construction presented in \cite{Shioda-Inose} which produces a K3 surface as quotient of the product surface $E_i\times E_i\simeq E_i\times \Cfu$.\\
Let $E_i$ be the curve $v^2=u(u^2+1)$ and $\alpha_E$ be the automorphism $\alpha_E:(u,v)\mapsto (-u,iv)$.\\
A generator of $H^{1,0}(E_i)$ is $\eta_E:=du/v$ and we observe that $\alpha_E^*(\eta_E)=i\eta_E$.\\
The surface $(E_i\times E_i)/\alpha_E^3\times \alpha_E$ is a singular surface and its desingularization is a K3 surface, described in   \cite[Lemma 5.2]{Shioda-Inose}. Here we give explicitly the quotient map.\\
Let the coordinates of the first copy of $E_i$ be $(u_1,v_1)$ and the coordinates of the second copy of $E_i$ be $(u_2,v_2)$.
Let  us consider the surface defined by the equation \begin{equation}\label{equation: Si surface}y^2=x^3+t^3(t+1)^2x.\end{equation}
We observe that $u_1v_2^2u_2^2$, $v_1v_2^3u_2^3$, $u_2^2$ are invariant under the action of $\alpha_E^3\times \alpha_E$. The quotient map $\pi:E_i\times E_i\ra (E_i\times E_i)/(\alpha_E^3\times \alpha_E)$ is
$$\pi:((u_1,v_1);(u_2,v_2))\mapsto (x:=u_2v_1^2u_1^2,y:=v_2v_1^3u_1^3,t:=u_1^2),$$
indeed it is clear that $x,y,t$ satisfy the equation \eqref{equation: Si surface} and if two points of $E_i\times E_i$ $((u_1,v_1);(u_2,v_2))\neq((u'_1,v'_1);(u'_2,v'_2))$
are such that $\pi((u_1,v_1);(u_2,v_2))=\pi((u_1,v_1);(u_2,v_2))$, then there are the following possibilities:
$$
\begin{array}{llll}u_1=u'_1,&v_1=-v'_1,&u_2=u'_2,&v_2=-v'_2\\
u_1=-u'_1,&v_1=iv'_1,&u_2=-u'_2,&v_2=-iv'_2\\
u_1=-u'_1,&v_1=-iv'_1,&u_2=-u'_2,,&v_2=iv'_2
\end{array}$$
thus $\pi$ identifies exactly the points in the orbit of $\alpha_E^3\times \alpha_E$ on $E_i\times E_i$.\\
The elliptic fibration \eqref{equation: Si surface} is the one described in \cite{Shioda-Inose}. The desingularization of such an elliptic fibration has 2 fibers of type $III^*$ and one of type $I_0^*$. The K3 surface $X_{SI}$, desingularization of the surface defined by the equation \eqref{equation: Si surface}, has Picard number 20 and hence its moduli space is zero dimensional.\\
The surface with equation \eqref{equation: Si surface} admits an automorphism $(x,y,t)\ra (-x,iy,t)$ (induced by $\alpha_E\times 1_C$), which extends to a non symplectic automorphism $\alpha_{X}$ of $X_{SI}$ such that $\alpha_X(\eta_{X})=i\eta_X$, where $\eta_X$ is a generator of $H^{2,0}(X_{SI})$.
The transcendental lattice of $X_{SI}$ is isometric to the lattice $\langle 2\rangle^2$ (cf.\ \cite{Shioda-Inose}).

\subsection{The cases $g=2,3$}\label{construction K3}
In the following we will denote by $H^{a,b}(\Cf)_{\lambda}$ the eigenspace of $H^{a,b}(\Cf)$ for the eigenvalue $\lambda$ with respect to the action of $\alpha_C$ on $H^{a,b}(\Cf)$ and by $H^{a,b}(\Sf)_{\lambda}$ the eigenspace of $H^{a,b}(\Sf)$ for the eigenvalue $\lambda$ with respect to the action of $\alpha_S$ (cf.\ Definition \ref{defi: alphaS}) on $H^{a,b}(\Sf)$.\\

Let us consider the product surface $E_i\times \Cf$ and the
automorphism $\alpha_E\times \alpha_C$ on it. One has:
\begin{eqnarray}\label{formula: H20, H11 of the quotient of E Cf }\begin{array}{ll}\left(H^{2,0}(E_i\times
\Cf)\right)^{\alpha_E\times \alpha_C}=&\left(H^{1,0}(E_i)\otimes
H^{1,0}(\Cf)\right)^{\alpha_E\times\alpha_C}=\langle
\eta_E\wedge \omega_C\rangle,\\ \\ \left(H^{1,1}(E_i\times
\Cf)\right)^{\alpha_E\times \alpha_C}=&\left\{\begin{array}{ll}
 \begin{array}{l}\left(H^{0,0}(E_i)\otimes H^{1,1}(\Cf)\right)\oplus\left(H^{1,1}(E_i)\otimes H^{0,0}(\Cf)\right)\\
  \oplus \langle\eta_E\wedge \overline{\eta_C}\rangle \oplus \langle\overline{\eta_E}\wedge \eta_C\rangle\end{array}&\mbox{ if }g=2,\\ \\
 \begin{array}{l}\left(H^{0,0}(E_i)\otimes H^{1,1}(\Cf)\right)\oplus\left(H^{1,1}(E_i)\otimes H^{0,0}(\Cf)\right)\\ \oplus 
 \langle \eta_E\wedge\overline{\eta_C} \rangle \oplus\langle \eta_E\wedge\overline{\nu_C}\rangle\oplus \langle
 \overline{\eta_E}\wedge\eta_C \rangle \oplus \langle \overline{\eta_E}\wedge\nu_c\rangle\end{array}&\mbox{ if
 }g=3,\end{array}\right.\end{array}\end{eqnarray}
where $\overline{\eta_C}$ is the conjugate of $\eta_C$, etc.
In particular $$dim (H^{1,1}(E_i\times \Cf)^{\alpha_E\times
\alpha_C})=\left\{\begin{array}{ll}4&\mbox{ if }g=2\\6&\mbox{ if
}g=3\end{array}\right.. $$

Let us consider the elliptic fibration \begin{equation}\label{equation: mathcal E}\mathcal{E}:\ \
y^2=x^3+xsf_g(s)^2.\end{equation} The map $\pi:E_i\times \Cf\ra \mathcal{E}$
defined by
$$\pi:((u,v);(z,r))\mapsto (x:=uz^2, y:=vz^3,s:=r^2)$$
is the quotient map $E_i\times \Cf\ra \left(E_i\times \Cf\right)/\alpha_E\times \alpha_C$ (this can be proved as in Section \ref{section: Shioda Inose}).\\
By standard results on the singularities of an elliptic fibration we easily deduce the kind of singularities of $\mathcal{E}$ by the zeroes of $sf_g(s)^2$ (cf.\ \cite[Table IV.3.1]{miranda elliptic pisa}): we have one singularity of type $A_1$ over $s=0$, $g$ singularities of type $D_4$ over the zeroes of $f_g(s)$ and one singularity over $s=\infty$ which is either of type $E_7$ if $g=2$ or of type $A_1$ if $g=3$. These singularities are A-D-E singularities, hence there exists a desingularization of $\left(E_i\times \Cf\right)/\left(\alpha_E\times \alpha_C\right)$ which is a K3 surface.
\begin{defi}\label{defi: alphaS} Let $\Sf$ be the desingularization of $\left(E_i\times \Cf\right)/\left(\alpha_E\times \alpha_C\right)$ and let $\alpha_S$ be the automorphism induced on $\Sf$ by the automorphism $\alpha_E$ on $E_i$.
\end{defi}
The automorphism $\alpha_S$ is $\alpha_S:(x,y,s)\mapsto (-x,iy,s)$. Since $H^{2,0}(\Sf)=\C \eta_S$, with $\eta_S=\eta_E\wedge \omega_C$, $\alpha_S(\eta_S)$ acts on $H^{2,0}(\Sf)$ as a multiplication by $i$.

\begin{rem}{\rm In case $g=1,2$ one can also consider the quotient $R_{f_g}:=\widetilde{(E_i\times \Cf)/(\alpha_E\times \alpha_C^3)}$ which admits an equation of the form $y^2=x^3+xs^3f_g(s)^2$. The quotient map generalizes the one considered in the previous section and is $$((u;v);(z;r))\ra (x:=uz^2r^2;y:=vz^3r^3,s:=r^2).$$ In case $g=1$, $R_{f_1}=X_{SI}$, in case $g=2$, the desingularization of $R_{f_2}$ is isomorphic to $\Sfd$
(it suffices to perform on $\mathbb{P}^1_s$ the automorphism $s\mapsto 1/s$). }\end{rem}

\begin{prop}\label{prop: elliptic fibration on E}
For a generic choice of the polynomial $f_g$, the K3 surface $\Sf$ admits an elliptic fibration with one fiber of type $III$ (two tangent rational curves) over $s=0$, $g$ fibers of type $I_0^*$ (with dual graph $\widetilde{D_4}$) over the zeroes of $f_g(s)$ and one reducible fiber over infinity which is either of type $III^*$ (with dual graph $\widetilde{E_7}$) if $g=2$, or of type $III$ if $g=3$. This fibration has no sections of infinite order, has a section of order 2, $\sigma$, and a rational curve $\beta$ as bisection.\end{prop}
\bprf The equation \eqref{equation: mathcal E} of $\mathcal{E}$ is the equation of the quotient surface $\left(E_i\times \Cf\right)/\left(\alpha_E\times \alpha_C\right)$. The desingularization of $\mathcal{E}$ gives an elliptic fibration with the singular fibers described (cf.\ \cite[Table IV.3.1]{miranda elliptic pisa}).\\ The equation of $\mathcal{E}$ can be written as $wy^2=x^3+xw^2sf_g(s)^2$. With these coordinates the zero section is $O:s\ra (x(s)$:$y(s)$:$w(s))=(0$:$1$:$0)$ and $\Sf$ admits a section of order 2, $\sigma:s\mapsto (0$:$0$:$1)$. Moreover it is clear that the fibration on  $\Sf$ admits a bisection $\beta$ with equation $x^2=-sf_g(s)^2$. This is a 2:1  cover of $\mathbb{P}^1$ branched over $s=0$ and $s=\infty$,  hence $\beta$ is a rational curve.\\
We have now to prove that there are no sections with infinite order.\\
The trivial lattice (i.e.\ the lattice generated by the class of the section $O$ and by the classes of the irreducible components of reducible fibers) of this elliptic fibration on $\Sf$ has rank equal to 18 if $g=2$ and to 16 if $g=3$. Hence the rank of the transcendental lattice $T_{\Sf}$ of $\Sf$ is \begin{eqnarray}\label{formula: rank Tg<}rank(T_{\Sf})\leq\left\{\begin{array}{l} 4\ \mbox{ if }g=2,\\ 6\ \mbox{ if }g=3\end{array}\right..\end{eqnarray} The family of curves $\Cf$ has $g-1$ moduli ($f_g$ has $g$ zeroes, and one of them can be chosen to be 1, so there are $g-1$ free zeroes of $f_g$) and since $E_i$ is rigid, the surface $\Sf$ has $g-1$ moduli. All the members of the family admit $\alpha_S$ as non symplectic automorphisms of order 4. The moduli of a family of K3 surfaces admitting a non symplectic automorphism of order $4$ is rank$(T)/2-1$, where $T$ is the transcendental lattice. Since the family $\Sf$ is at least a subfamily of the family of K3 surfaces admitting a non symplectic automorphism of order 4, the number of moduli of $\Sf$ is at most $rank(T_{\Sf})/2-1$, hence
\begin{eqnarray}\label{formula: rank Tg>}\frac{rank(T_{\Sf})}{2}-1\geq \left\{\begin{array}{l}1\ \mbox{ if }g=2\\ 2\ \mbox{ if } g=3\end{array}\right..\end{eqnarray} Comparing \eqref{formula: rank Tg<} and \eqref{formula: rank Tg>} one obtains that $rank(T_{\Sf})=g-1$. In particular this implies that the rank of the N\'eron Severi group of $\Sf$ is the same as the rank of the trivial lattice of the elliptic fibration on $\Sf$ and hence there are no sections of infinite order.\eprf 
\begin{rem}\label{rem: T of Sfg}{\rm Since the classes of the irreducible components of the reducible fibers and of the sections of an elliptic fibration generate the N\'eron Severi group of the surface, one can directly compute the discriminant form of the N\'eron--Severi group of each of these K3 surfaces. By this computation one deduces that the transcendental lattice of $S_{f_3}$ is isometric to $U(2)^2\oplus \langle -2\rangle^2$ and the transcendental lattice of $S_{f_2}$ is isometric to $U(2)\oplus \langle 2\rangle\oplus\langle -2\rangle$.}\end{rem}

\subsection{The order 4 non symplectic automorphism on $\Sf$.}\label{non symplectic 4} 
\begin{prop}\label{prop: non sympl autom on S} In the case $g=2$, the fixed locus of $\alpha_S$ is made up of 3 rational curves and 10 points and the fixed locus of $\alpha_S^2$ is made up of 8 rational curves (3 of them are fixed by $\alpha_S$, the others are invariant for $\alpha_S$).\\
In the case $g=3$, the fixed locus of $\alpha_S$ is made up of 2 rational curves and 8 points and the fixed locus of $\alpha_S^2$ is made up of 6 rational curves (2 of them are fixed by $\alpha_S$, the others are invariant for $\alpha_S$). 
\end{prop}
\bprf {\it The fixed locus of $\alpha_S^2$.} The automorphism $\alpha_S^2$ is a non symplectic involution on a K3 surface, hence its fixed locus is either empty or made up of disjoint curves (cf.\ \cite{nikulin order 2}). The sections $O$ and $\sigma$ and the bisection $\beta$ are fixed rational curves.\\
We now describe the action of $\alpha_S^2$ on the irreducible components of the reducible fibers. We will call $C_0^i$ the component of the $i$-th reducible fiber which meets the zero section. For the fibers of type $III$ we call $C_1^i$  the other component.\\
The 2-torsion section meets the fibers of type $III$ in the component $C_1$. The component $C_0$ (resp.\ $C_1$) meets the zero section (resp.\ the section $\sigma$), which is a fixed curve and hence cannot be fixed (the fixed locus consists of disjoint curves). The point of intersection of $C_0$ and $C_1$ is fixed, but since the fixed locus does not contain isolated fixed points, this point is on the bisection $\beta$ (one can also directly check that the bisection $\beta$ meets the two rational curves of the fiber in their intersection point).\\
The fibers of type $I_0^*$ are made up of five rational curves with the following configuration:
$$
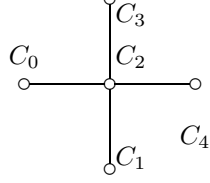
\DynkinDD{C_0}{C_1}{C_2}{C_3}{C_4}
$$
The section $O$ meets the component $C_0$, the section $\sigma$ the component $C_1$ and the bisection $\beta$ the components $C_3$ and $C_4$. Since the fixed locus is smooth, one obtains that the components $C_i$,
$i=0,1,3,4$ are invariant, but they are not pointwise fixed. The points of intersection of $C_2$ with the curves $C_i$,
$i=0,1,3,4$ are fixed and hence $C_2$ is pointwise fixed.\\
The fiber of type $III^*$ is made up of eight rational curves
with the following configuration:
$$
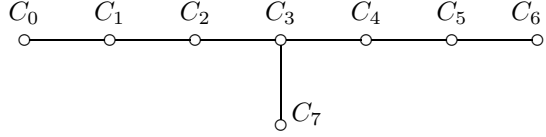
\DynkinEIII{C_0}{C_1}{C_2}{C_3}{C_4}{C_5}{C_6}{C_7}
$$
The section $O$ meets $C_0$, $\sigma$ meets $C_6$ and the bisection $\beta$ is tangent to $C_7$. Arguing as before one concludes that the components $C_1$, $C_3$, $C_5$ are pointwise fixed, the others are invariant but not fixed.\\
So the curves fixed by $\alpha_S^2$ are: $O$, $\sigma$, $\beta$, the components $C_2$ on each fibers of type $I_0^*$ and, if $g=2$, the components $C_j$, $j=1,3,5$ on the fiber of type $III^*$.\\
{\it Fixed locus of $\alpha_S$.} The sections $O$ and $\sigma$ are fixed by $\alpha_S$. The bisection $\beta$ is invariant, but not fixed. Since it is a rational curve, the restriction of $\alpha_S$ on it has two fixed points: these are the points of intersection of $\beta$ with the reducible fibers over $s=0$ and $s=\infty$ (i.e.\ on the fibers which are not of type $I_0^*$). Thus the two points of intersection of $\beta$ with each fibers of type $I_0^*$ are switched by $\alpha_S$ and hence on the fiber of type $I_0^*$ the components $C_3$ and $C_4$ are switched, so $C_2$ can not be pointwise fixed but the intersection points between $C_2$ and the components $C_0$ and $C_1$ are fixed.\\
On the fiber of type $III^*$ the point of intersection of $C_7$ and $\beta$ is a fixed point. The component $C_3$ is pointwise fixed and the components $C_1$ and $C_5$ are invariant. It is possible that
they are either pointwise fixed or invariant but not pointwise fixed. To understand which case occurs we consider the Lefschetz fixed point formula:
$$\chi(Fix_{\alpha_S}(S))=\sum_r(-1)^r\mbox{trace}({\alpha_S}_{|H^r(S,\C)}),$$
where $Fix_{\alpha_S}(S)$ is the fixed locus of $\alpha_S$ on $S$ and $\chi$ is the topological Euler characteristic.\\
Here we are in case $g=2$ (because for $g=3$ there is no fiber of type $III^*$). From the description of $\alpha_S$ it is clear that $\alpha_S$ acts with eigenvalue $-1$ on the classes $C_3^{(2)}-C_4^{(2)}$ and $C_3^{(3)}-C_4^{(3)}$ (assuming that the second and the third fibers are the fibers of type $I_0^*$) and acts trivially on the other classes in the N\'eron--Severi group.
Hence the eigenspace for the eigenvalue $-1$ on $NS(\Sfd)\otimes\C$ has dimension 2 and the eigenspace for the eigenvalue $1$ has dimension $16$. Moreover $\alpha_S$ acts as the multiplication by $i$ on the period $\eta_{\Sfd}$, so it acts with eigenvalue $i$ on a 2-dimensional subspace of $T_S\otimes \C$ and with eigenvalue $-i$ on the orthogonal 2-dimensional subspaces of $T_S\otimes \C$ (cf.\ \cite{Nikulin symplectic}). The action of $\alpha_S$ on $H^0(\Sfd,\C)$ and on $H^4(\Sfd,\C)$ is trivial and $H^1(\Sfd \C)=H^3(\Sfd,\C)=0$ since $\Sfd$ is a K3 surface. Hence the right side in the Lefschetz fixed point formula is $1+16+2(-1)+2i+2(-i)+1=16.$ We recall that the Euler characteristic of a point is 1 and of a rational curve is 2. By the description of $\alpha_S$ either $Fix_{\alpha_S}(\Sfd)=\left(\coprod_{i=1}^6 p_i\right)\coprod\left(\coprod_{j=1}^5C_j\right)$ or $Fix_{\alpha_S}(\Sfd)=\left(\coprod_{i=1}^{10} p_i\right)\coprod\left(\coprod_{j=1}^3 C_j\right)$ where $p_i$ are points and $C_j$ are rational curves. It follows that $Fix_{\alpha_S}(\Sfd)=\left(\coprod_{i=1}^{10} p_i\right)\coprod\left(\coprod_{j=1}^3 C_j\right)$ and hence the components $C_1$ and $C_5$ on the fiber of type $III^*$ are not pointwise fixed.\\
To recap, if $g=2$ the fixed locus of $\alpha_S$ is made up of 3 rational curves ($O$, $\sigma$, and the component $C_3$ on the fiber of type $III^*$), and 10 points (1 on the fiber of type $III$, 2 on each fiber of type $I_0^*$ and 5 on the fiber of type $III^*$); if $g=3$ the fixed locus of $\alpha_S$ is made up of 2 rational curves ($O$, $\sigma$), and 8 points (1 on each fiber
of type $III$ and 2 on each fiber of type $I_0^*$).\eprf
$ $\\
In \eqref{formula: H20, H11 of the quotient of E Cf } we computed
the invariant part of $H^{2,0}(E_i\times\Cf)$ and $H^{1,1}(E_i\times \Cf)$ under $\alpha_E\times \alpha_C$. These spaces pass to the quotient. In the previous proposition we calculated explicitly the action of $\alpha_S$ on the curves $C_i^j$, which are the curves arising from the resolution of the singular quotient $\left(E_i\times\Cf\right)/\left(\alpha_E\times \alpha_C\right)$. These curves are classes in $NS(\Sf)\subset H^{1,1}(\Sf)$. Hence here we can give the complete decomposition in eigenspaces of the second cohomology group of $H^2(\Sf,\C)$:

\begin{align}\label{formula: decomp. eigensp. H2(S)}
&H^{2}(\Sf,\C)_i&=&\left\{\begin{array}{lll}\langle \eta_E\wedge\omega_C
\rangle\oplus
\langle\eta_E\wedge\overline{\eta_C}\rangle&\simeq
H^{2,0}(\Sf)\oplus\left(H^{1,0}(E_i) \otimes H^{0,1}(\Cf)_{-i}
\right)&\mbox{ if }g=2\\\langle\eta_E \wedge\omega_C
\rangle\oplus
\langle\eta_E\wedge\overline{\eta_C}\rangle\oplus
\langle\eta_E\wedge\overline{\nu_C}\rangle&\simeq
H^{2,0}(\Sf)\oplus \left(H^{1,0}(E_i)\otimes H^{0,1}(\Cf)_{-i}
\right)&\mbox{ if
}g=3\end{array}\right.\\
\nonumber
&H^{2}(\Sf,\C)_{-1}&=&H^{1,1}(\Sf)_{-1}=\left\{\begin{array}{lll}\langle
C_3^{(2)}-C_4^{(2)},C_3^{(3)}-C_4^{(3)} \rangle\otimes\C&\mbox{ if
}g=2\\\langle
C_3^{(2)}-C_4^{(2)},C_3^{(3)}-C_4^{(3)},C_3^{(4)}-C_4^{(4)},
\rangle\otimes\C&\mbox{ if
}g=3\end{array}\right.\\
\nonumber
&H^{2}(\Sf,\C)_{-i}&=&\left\{\begin{array}{lll}\langle
\overline{\eta_E}\wedge \overline{\omega_C}\rangle\oplus
\langle\overline{\eta_E}\wedge\eta_C\rangle&\simeq
H^{0,2}(\Sf)\oplus\left(H^{0,1}(E_i)\otimes H^{1,0}(\Cf)_{-i}
\right)&\mbox{ if }g=2\\\langle  \overline{\eta_E}
\wedge\overline{\omega_C}\rangle\oplus
\langle\overline{\eta_E}\wedge\eta_C\rangle\oplus
\langle\overline{\eta_E}\wedge\nu_C\rangle&\simeq
H^{0,2}(\Sf)\oplus \left(H^{0,1}(E_i)\otimes H^{1,0}(\Cf)_{-i}
\right)&\mbox{ if
}g=3\end{array}\right.\\
\nonumber
&H^{2}(\Sf,\C)_{1}&=&H^{1,1}(\Sf)_{1}=\left\{\begin{array}{lll}\langle
C_3^{(2)}-C_4^{(2)},C_3^{(3)}-C_4^{(3)} \rangle
^{\perp_{NS(\Sf)}}\otimes\C&\mbox{ if }g=2\\\langle
C_3^{(2)}-C_4^{(2)},C_3^{(3)}-C_4^{(3)},C_3^{(4)}-C_4^{(4)},
\rangle^{\perp_{NS(\Sf)}}\otimes\C&\mbox{ if
}g=3\end{array}\right.
\end{align}
where, if $L$ is a sublattice of $M$, we indicate with $L^{\perp_M}$ the sublattice of $M$ which is orthogonal to $L$.

\begin{rem}{\rm The previous computation can also be applied to the surface $X_{SI}$. By the description of the action of the non symplectic automorphism $\alpha_S$ on the singular fibers of the fibrations on $\Sf$ one deduces the properties of the fixed locus of automorphism $\alpha_{X}$: $\alpha_X$ fixes 4 rational curves (one for each fiber of type $III^*$, the zero sections and the 2-torsion section) and 12 points (5 on each fibers of type $III^*$ and 2 on the fiber of type $I_0^*$), and $\alpha_X^2$ fixes 10 rational curves (4 of them are fixed by $\alpha_X$, the other are invariant for $\alpha_X$). Moreover one can describe explicitly the eigenspaces for the action of $\alpha_X$ on $H^2(X_{SI},\Z)$ as done for the surfaces $\Sf$.\erem} \end{rem}

\section{Isotrivial fibrations with section and the family $Q_{b^2}$ of K3 surfaces}\label{subsection: isotrivial} In the previous section we described K3 surfaces admitting an isotrivial fibration with section such that the smooth fibers are isomorphic to $E_i$. Each such a fibration has a minimial Weierstrass equation of type $y^2=x^3+xa(t,s)$ where $a(t,s)$ is a homogeneous polynomial of degree 8 and the maximal multiplicity of its zeroes is 3 (the conditions on the degree and on the multiplicity of the zeroes guarantee that the elliptic fibration is a K3 surface, \cite{miranda elliptic pisa}). 
\begin{defi}
Let $Q_a$ be the K3 surface defined by $y^2=x^3+xa(t,s)$  and $\alpha_Q$ be the automorphism induced by the automorphism $\alpha_E$ acting on the fibers, i.e.\ $\alpha_Q:(x,y,t,s)\mapsto (-x,iy,t,s)$.\end{defi}
In order to construct families of Calabi--Yau 3-folds without maximal unipotent monodromy, we require that $\alpha_Q$ and its powers do not fix curves of positive genus. As in previous section one observes that the fixed loci of $\alpha_Q$ and $\alpha_Q^2$ are made up of points, components of the reducible fibers (which are always rational curves), the zero section, the torsion section $\sigma$:$(t$:$s)\mapsto (0,0,(t$:$s))$ (these are rational curves) and the bisection $x^2=a(t,s)$. Hence $\alpha_Q$ and $\alpha_Q^2$ fix only points and rational curves if and only if $x^2=a(t,s)$ is not a curve of positive genus. If $a(t,s)$ is not a square, then $x^2=a(t,s)$ is a hyperelliptic curve and its genus depends on the zeroes with an odd multiplicity of the polynomial $a(t,s)$. Recalling that the maximum degree of a zero in $a(t,s)$ is 3 one sees that there are only 3 possibilities for $a(t,s)$, which correspond to the surfaces $X_{SI}$, $S_{f_2}$ and $S_{f_3}$ constructed in previous sections.\\
The only other possibility is that $a(t,s)=b^2(t,s)$ where $b(t,s)$ is a homogeneous polynomial of degree 4 without multiple roots. In this case the bisection $x^2=a(t,s)$ splits in the union of two sections, which are rational curves since the base curve of the fibration is $\mathbb{P}^1$. Thus we have another family $Q_{b^2}$ of K3 surfaces. 
\begin{prop}\label{prop: Qb2} Let $E$ be the elliptic curve $E: v^2=b(t,s)$, $deg(b)=4$, and let $\iota_E$ be the automorphism $\iota_E:(v,t,s)\mapsto (-v,t,s)$. The desingularization of the product $(E_i\times E)/(\alpha_E^2\times \iota_E)$ is the K3 surface $Q_{b^2}$ admitting the isotrivial elliptic fibration $y^2=x^3+xb^2(t,s)$.\\
This fibration  has four fibers of type $I_0^*$ and three 2-torsion sections.\\ 
The automorphism $\alpha_Q$ fixes 2 rational curves and 8 points, the automorphism $\alpha_Q^2$ fixes 8 rational curves (2 of them are fixed by $\alpha_Q$, 4 are invariant for $\alpha_Q$ and 2 are switched by $\alpha_Q$). \end{prop}
\bprf The proof is similar to the ones of Propositions \ref{prop: elliptic fibration on E} and \ref{prop: non sympl autom on S}.\eprf
\begin{rem}\label{rem: T of Qb}{\rm The transcendental lattice of $Q_{b^2}$ is isometric to $U(2)\oplus U(2)$ (cf.\ Remark \ref{rem: T of Sfg}).}\end{rem}

\section{The $4:1$ covers of $\mathbb{P}^2$ and the family $R$ of K3 surfaces}\label{section: geometrical description}
We discuss a different projective model for the K3 surfaces $\Sf$ (Proposition \ref{prop: 4:1 cover known}) and give another family $R$ of K3 surfaces (Proposition \ref{prop: a 3 dimensional family}) which admits a similar model.\\
In \cite{michela Z4 action} a one dimensional family of K3 surfaces $Z$ admitting a non symplectic automorphism is constructed as 4:1 cover of $\mathbb{P}^2$ branched along a quartic which is the intersection of a conic and two lines. In \cite[Proposition 3.1]{michela Z4 action} it is proved that this family admits an isotrivial fibration with general fiber isomorphic to $E_i$ and the equation of the curve $B$ which trivializes the fibration is given. In other words the curve $B$ described in \cite{michela Z4 action} is a curve such that there exists a finite group $G$ such that the desingularization of $(B\times E_i)/G$ is the K3 surface $Z$.  From the proof of \cite[Proposition 3.1]{michela Z4 action} it is evident that the curve $B$ is the curve $\Cfd$ and the group $G$ is $\Z/4\Z$ which acts on $E_i\times \Cfd$ as described in the Section \ref{construction K3}. So the one dimesional family of surfaces which are 4:1 cover of $\mathbb{P}^2$ branched along a conic and two lines is the family $\Sfd$. Moreover, the non symplectic automorphism described in the Section \ref{non symplectic 4} is the automorphism of order 4 associated to the 4:1 cover of $\mathbb{P}^2$. We observe that the automorphism associated to the cover has to fix the three rational curves which are in the branch locus. Moreover it fixes 10 points. Indeed, to obtain a smooth model of the 4:1 cover of $\mathbb{P}^2$, one has to blow up the singular points of the branch quartic, introducing three rational curves with an $A_3$ configuration for each singular points. The intersection points between the curves of each $A_3$ configuration are fixed by the automorphism. Since the branch quartic has 5 singular points, the automorphism has 10 fixed points. The square of this automorphism fixes the rational curves fixed by the automorphism, and the central components of the curves in the $A_3$ configuration.\\
Here we prove that also the surfaces ${X}_{SI}$ and $\Sft$ admit a similar description.
\begin{prop}\label{prop: 4:1 cover known} The 2-dimensional family of K3 surfaces which are 4:1 covers of $\mathbb{P}^2$ branched along the union of two conics in general position is the family $\Sft$. The K3 surface which is the 4:1 cover of $\mathbb{P}^2$ branched along the union of four lines in general position is the surface $X_{SI}$.\end{prop}
\bprf Up to a choice of the coordinates of $\mathbb{P}^2$ one can assume that the intersection between two conics are the points $(1:0:0)$, $(0:1:0)$, $(0:0:1)$, $(1:1:1)$. So one can always assume that the K3 surface which is the $4:1$ cover of $\mathbb{P}^2$ branched along two conics is the 2-dimensional family
$$s^4=(x_1x_2+ax_2x_3+(a-1)x_1x_3)(x_1x_2+bx_2x_3+(b-1)x_1x_3).$$
If we now consider the pencil of 
lines through $(0:0:1)$, i.e.\ $x_2=\lambda x_1$, this gives a pencil of elliptic curves on the surface with $j$-invariant equal to 1728 (i.e.\ an isotrivial fibration on the surface with generic fiber isomorphic to $E_i$). To trivialize this fibration we consider the intersection of the pencil with the branch quartic and we put $x_1=1$, so the equation of the surface becomes
$$s^4=(\lambda+a\lambda x_3+(a-1)x_3)(\lambda+b\lambda x_3+(b-1)x_3).$$
Now we have that $$(x_3-h)(x_3-k)=\frac{1}{4}(h-k)^2(z^2-1),\mbox{ where }z\mbox{ is such that } x_3=\frac{1}{2}(h-k)z+\frac{1}{2}(h+k),$$
so putting $h=\frac{\lambda}{a\lambda+a-1}$ and $k=\frac{\lambda}{b\lambda+b-1}$ one obtains 
$$s^4=(a\lambda+a-1)(b\lambda+b-1)\frac{\lambda^2(\lambda+1)^2(b-a)^2}{4(a\lambda+a-1)^2(b\lambda+b-1)^2}(z^2-1)=\frac{\lambda^2(\lambda+1)^2(b-a)^2}{4(a\lambda+a-1)(b\lambda+b-1)}(z^2-1).$$
Considering the base change $B_{(\rho,\lambda)}\mapsto \mathbb{P}^2_\lambda$ with $$B_{(\rho,\lambda)}:= \left\{\rho^4=\frac{\lambda^2(\lambda+1)^2(b-a)^2}{4(a\lambda+a-1)(b\lambda+b-1)}\right\}$$ and the change of coordinates $\tau:=\frac{\lambda(\lambda+1)(b-a)^{1/2}}{\rho(4ab)^{1/4}}$ one obtains that the trivializing curve $B_{(\tau, \lambda)}$ is $$\tau^4=\lambda^2(\lambda+1)^2(\lambda-a')(\lambda-b'),\ \mbox{ where }a':=\frac{a}{a-1},\ \ b':=\frac{b}{b-1}.$$ This is the curve $\Cft$ with the equation given in \eqref{equation D} (assuming that one zero of $l_2$ is at infinity). So the surfaces which are $4:1$ cover of $\mathbb{P}^2$ branched along a conic and two lines are the quotient of $E_i\times \Cft$ by the automorphism $((s,z),(\tau, \lambda))\mapsto((is,z),(-i\tau,\lambda))$, hence the family of K3 surfaces which admits such a 4:1 cover of $\mathbb{P}^2$ is $\Sft$.\\
The proof in the case the branch locus is made up of four lines is very similar. We only observe that the equation of the surface is $s^4=x_1x_3(x_1-x_3)(x_2-x_3)$.
\eprf

\begin{prop}\label{prop: a 3 dimensional family} There exists a 3-dimensional family of K3 surfaces $R$ which are $4:1$ covers of $\mathbb{P}^2$ branched along a quartic curve with three nodes. All the K3 surfaces in this family admit an isotrivial fibration with generic fiber isomorphic to $E_i$ and the order four automorphism, $\alpha_R$,  induced on $R$ by $\alpha_E$ is the cover automorphism. The automorphism $\alpha_R$ fixes 1 rational curve (the branch curve) and 6 points, the automorphism $\alpha_Q^2$ fixes 4 rational curves (1 of them is fixed by $\alpha_R$, the others are invaraint for $\alpha_R$).\end{prop}
\bprf It is well known that the $4:1$ cover of $\mathbb{P}^2$ branched over a quartic curve with at most A-D-E singularities admits a desingularization which is a K3 surface. It is immediate to check that, up to an automorphism of $\mathbb{P}^2$, the irreducible plane quartic curves with three nodes are $V(x_1^2x_2^2+x_1^2x_3^2+x_2^2x_3^2+ax_1^2x_2x_3+bx_1x_2^2x_3+cx_1x_2x_3^2)$. Thus the family of $R$ is 3-dimensional. The equation of the $4:1$ cover is $s^4=x_1^2x_2^2+x_1^2x_3^2+x_2^2x_3^2+ax_1^2x_2x_3+bx_1x_2^2x_3+cx_1x_2x_3^2$ and the cover automorphism is $\alpha_R:(s,x_1,x_2,x_3)\mapsto (is,x_1,x_2,x_3,x_4)$. The generic line of the pencil $x_2=\lambda x_1$ intersects the quartic in the node $(0:0:1)$ and in two other points. This pencil of lines exhibits an elliptic fibration on $R$ with general fiber isomorphic to an elliptic curve with $j$ invariant equal to 1728 (we put $x_2=1$): 
\begin{equation}\label{eq: isotrivial without section} s^4=\left(1+c\lambda+c\lambda^2\right)\left(x_3^2-\lambda^2\left(\left(\frac{(a+b\lambda)}{2(1+c\lambda+c\lambda^2)}\right)^2+1\right)\right).\end{equation} 
For almost every values of $\lambda$ the Equation \eqref{eq: isotrivial without section} gives an elliptic curve isomorphic to $E_i$ with a non symplectic automorphism of order 4 $(s,x_3)\ra (is,x_3)$, so $\alpha_E$ induces exactly $\alpha_R$. The fixed locus can be computed as in case $S_{f_2}$.
\eprf
In particular there exists a curve $B$ and a finite group $G$ of automorphism of $B\times E_i$ such that $R$ is a desingularization of $\left(B\times E_i\right)/G$. 

The family of K3 surfaces admitting a non symplectic involution fixing four rational curves is a 6-dimensional family and the transcendental lattice of a generic member of this family is isomorphic to $U(2)^2\oplus\langle -2\rangle^4$  (cf.\ \cite{nikulin order 2}). It is the eigenspace $T$ of the eigenvalue $-1$ for the actions of the involution on $H^2$. If the non symplectic involution is the square of a non symplectic automorphism of order 4, the space $T$ splits in the sum of the two 4-dimensional eigenspaces of the eigenvalues $i$ and $-i$ for the action of the automorphims of order 4 on $H^2$. The $H^{2,0}$ of a member of this family lives in one of these two eigenspaces and we can assume it lives in the eigenspace of the eigenvalue $i$. Hence the family of K3 surfaces admitting a non symplectic automorphism of order 4 whose square fixes four rational curves, has dimension $4-1=3$ and the transcendental lattice of a general member is isomorphic to $U(2)^2\oplus\langle -2\rangle^4$. Since $R$ is a 3-dimensional family of K3 surfaces admitting an automorphism  $\alpha_R$ of order 4 such that $\alpha_R^2$ fixes four rational curves, these two families coincide and the transcendental lattice of $R$, $T_R$, is isomorphic to $U(2)^2\oplus \langle -2\rangle^4$. By \cite{nikulin order 2}, $\alpha_R^2$ acts as the identity on the N\'eron--Severi group of the surface and as $-1$ on the transcendental lattice. Since $\alpha_R$ acts on $H^{2,0}(R)$ as a multiplication by $i$, the eigenspace $H^{2}(R,\C)_i$ of the eigenvalue $i$ for $\alpha_R$ decomposes in $H^{2,0}(R)\oplus H^{1,1}(R)_i$ and has dimension $4$, as said before. Similarly the eigenspace $H^{2}(R,\C)_{-i}$ has dimension 4 and decomposes in $H^{0,2}(R)\oplus H^{1,1}(R)_{-i}$. The sum of the eigenspaces with eigenvalues $+1$ and $-1$ is $NS(R)\otimes \C$. Let $x$ be the dimension of the invariant sublattice of $H^{1,1}(R)$ for $\alpha_R$, and so $14-x$ is the dimension of the anti-invariant sublattice of $H^{1,1}(R)$. Let $C$ be the rational curve and $p_i$ be the six points  fixed by $\alpha_R$. By the Lefschetz fixed point formula we obtain: $$8=\chi\left(\left(\coprod_{i=1}^6 p_i\right)\coprod\left(\coprod C\right)\right)=1+x+(14-x)(-1)+4i+4(-i)+1=-12+2x.$$ 
So $x=4$. Thus the dimensions of the eigenspaces of $H^2(R,\Q)$ for the action of $\alpha_R$ are: $$\begin{array}{ll}dim\left(H^{2}(R,\C)_{i}\right)=dim(H^{2,0}(R)\oplus H^{1,1}(R)_i)=4,&dim\left(H^{2}(R,\C)_{-1}\right)=dim(H^{1,1}(R)_{-1})=4,\\dim\left(H^{2}(R,\C)_{-i}\right)=dim(H^{2,0}(R)\oplus H^{1,1}(R)_{-i})=4,&dim\left(H^{2}(R,\C)_{1}\right)=dim(H^{1,1}(R)_1)=10.
\end{array}$$

\begin{rem}{\rm In \cite{nikulin order 2} it is proved that the families of K3 surfaces admitting a non symplectic involution fixing $k$ rational curves, $k=10,8,6,4$, is $(10-k)$-dimensional and the generic member of each of these families has transcendental lattice isometric to $U(2)$ if $k=10$, either to $U(2)^{\oplus 2}$ or to $U(2)\oplus\langle 2 \rangle\oplus \langle -2\rangle$ if $k=8$, to $U(2)^{\oplus 2}\oplus \langle -2\rangle^{\oplus 2}$ if $k=6$, to $U(2)^{\oplus 2}\oplus \langle -2\rangle^{4}$ if $k=4$. For each of these families we constructed a $((10-k)/2)$-dimensional subfamily admitting an automorphism of order 4 whose square is the involution fixing $k$ rational curves.}\end{rem}
\section{Families of Calabi Yau 3-folds.}\label{section: Calabi Yau} 
Once one has a family of K3 surfaces with a non symplectic automorphism of order 4 such that its fixed locus and the fixed locus of its square consist of rational curves and points, one can construct Calabi--Yau 3-folds as quotient of the product of these K3 surfaces with the elliptic curve $E_i$. In the previous section we constructed some families of K3 surfaces with the required automorphism. In Section \ref{subsec: Yf} we construct  families of Calabi--Yau 3-folds $\Yf$ considering the families of K3 surfaces $S_{f_g}$, giving all the details.  The construction in the other two cases (starting with the families of K3 surfaces $Q_{b^2}$ and $R$) is very similar and is presented in sections \ref{subsection: Wb} and \ref{subsection: 3 dim fam CY}.
\subsection{The desingularization $\Yf$ of the quotient $\left(E_i\times \Sf\right)/\left(\alpha_E^3\times \alpha_S\right)$}\label{subsec: Yf}
Let $\Vf$ be the (singular) 3-fold $\left(E_i\times \Sf\right)/\left(\alpha_E^3\times \alpha_S\right)$ and $\Yf$ its crepant desingularization.
\begin{rem}\label{rem: desingularization}{\rm The desingularization of $\Vf$ can be also viewed as the desingularization of the quotient $\left(E_i\times E_i\times \Cf\right)/\left(\langle\alpha_E^3\times \alpha_E\times 1,1\times \alpha_E\times\alpha_C\rangle\right)$ or of the quotient $\left(X_{SI}\times \Cf\right)/\left(\alpha_X\times \alpha_C\right)$. Indeed let
$S_{sing}:=\left(E_i\times \Cf\right)/\left(\alpha_E\times
\alpha_C\right)$ (resp.\ $X_{sing}:=\left(E_i\times E_i\right)/\left(\alpha_E\times
\alpha_E^3\right)$) and $\alpha_{S_{s}}$ (resp.\ $\alpha_{X_s}$) the automorphism induced by
$\alpha_E$. So we have birationals isomorphisms: 
{\small \begin{align}\nonumber\left(E_i\times
S_{sing}\right)/\left(\alpha_E^3\times \alpha_{S_{s}}\right)\simeq
\left(E_i\times E_i\times \Cf\right)/\left(\langle\alpha_E^3\times
\alpha_E\times 1,1\times \alpha_E\times\alpha_C\right\rangle)\simeq
\left(X_{sing}\times \Cf\right)/\left(\alpha_{X_s}\times
\alpha_C\right).\end{align}} Since the automorphisms
$\alpha_E^3\times \alpha_E\times 1$ and $1\times
\alpha_E\times\alpha_C$ commute on $E_i\times E_i\times \Cf$ we
can consider the quotient $\left(E_i\times E_i\times
\Cf\right)/\left(1\times \alpha_E\times \alpha_C\right)$ (resp.\ $\left(E_i\times E_i\times
\Cf\right)/\left(\alpha_E^3\times \alpha_E\times 1\right)$), and resolve
its singularities,  obtaining $E_i\times \Sf$ (resp.\ $X_{SI}\times \Cf $). Since the
automorphism $\alpha_E^3\times \alpha_E\times 1$ (resp.\ $1\times \alpha_E\times\alpha_C$) induces the
automorphism $\alpha_E\times \alpha_S$ (resp.\ $\alpha_X\times \alpha_C$) on $E_i\times \Sf$ (resp.\ on $X_{SI}\times \Cf$), the
desingularization $\Yf$ of the quotient $\Vf$ is a
desingularization of $\left(E_i\times E_i\times
\Cf\right)/\left(\alpha_E^3\times \alpha_E\times 1,1\times
\alpha_E\times\alpha_C\right)$ (resp.\ $X_{SI}\times \Cf/\left(\alpha_X\times \alpha_C\right)$).}\end{rem}
The holomorphic forms on the smooth part of the quotient $\Vf$ are the images of the $\alpha_E^3\times\alpha_S$-invariant holomorphic forms on $E_i\times \Sf$. Since $\alpha_S$ acts as the identity on $H^{0,0}(\Sf)$ and as the multiplication by $i$ on $H^{2,0}(\Sf)$ and $\alpha_E^3$ acts as the identity on $H^{0,0}(E_i)$ and as the multiplication by $-i$ on $H^{1,0}(E_i)$, we have: \begin{align}\label{formula: H30 Y}\begin{array}{ll}(H^{3,0}(E_i\times \Sf))^{\alpha_E^3\times\alpha_S}=&\left(H^{1,0}(E_i)\otimes H^{2,0}(\Sf)\right)^{\alpha_E^3\times \alpha_S}=\langle\eta_E\times \eta_S\rangle\\
(H^{2,0}(E_i\times \Sf))^{\alpha_E^3\times\alpha_S}=&\left(\left(
H^{0,0}(E_i)\otimes H^{2,0}(\Sf)\right)\oplus\left(
H^{1,0}(E_i)\otimes H^{1,0}(\Sf)\right)\right)^{\alpha_E^3\times \alpha_S}=0\\
(H^{1,0}(E_i\times \Sf))^{\alpha_E^3\times\alpha_S}=&\left(\left(
H^{0,0}(E_i)\otimes H^{1,0}(\Sf)\right)\oplus\left(
H^{1,0}(E_i)\otimes H^{0,0}(\Sf)\right)\right)^{\alpha_E^3\times
\alpha_S}=0.\end{array}\end{align} Moreover
\begin{align}\label{formula: H21 Y}\begin{array}{lll}H^{2,1}(E_i\times \Sf)^{\alpha_E^3\times \alpha_S}&=H^
{1,0}(E_i)\otimes \left(H^{1,1}(\Sf)\right)_i\\
H^{1,1}(E_i\times \Sf)^{\alpha_E^3\times \alpha_S}&=\left(H^{0,0}(E_i)\otimes  H^{1,1}(\Sf)_1\right)\oplus
\left(H^{1,1}(E_i)\otimes
H^{0,0}(\Sf)\right).\end{array}\end{align}
$$\mbox{So }dim(H^{2,1}(E_i\times
\Sf)^{\alpha_E^3\times\alpha_S})=\left\{\begin{array}{rr}1&\mbox{if
}g=2\\2&\mbox{if }g=3\end{array}\right.\ \ \mbox{ and
}dim(H^{1,1}(E_i\times
\Sf)^{\alpha_E^3\times\alpha_S})=\left\{\begin{array}{rr}17&\mbox{if
}g=2\\14&\mbox{if }g=3\end{array}\right..$$

The locus with a non trivial stabilizer for the action of $\alpha_E^3\times \alpha_S$ on $E_i\times \Sf$ is given by the product of the locus with a non trivial stabilizer for the action of $\alpha_S$ on $\Sf$ and of $\alpha_E^3$ on $E_i$. On $E_i$ there are two points fixed by $\alpha_E^3$, on which the local action is the multiplication by $-i$. The action of $\alpha_S$ can be linearized and diagonalized near the fixed locus (cf.\ \cite{Nikulin symplectic}), in particular near the fixed curves it is $diag(i,1)$ and near the isolated fixed points it is $diag(i^3,-1)$. Hence on $E_i\times \Sf$ the action of $\alpha_E^3\times\alpha_S$ linearizes as $diag(i,1,-i)$ near  the curves fixed by $\alpha_E^3\times\alpha_S$ and as $diag(-i,-1,-i)$ near the points fixed by $\alpha_E^3\times\alpha_S$.

\begin{prop}\label{prop: Hodge numbers CY} There exists a desingularization  $\Yf$ of $\left(E_i\times \Sf\right)/\left(\alpha_E^3\times \alpha_S\right)$ which is a smooth Calabi--Yau variety. Its Hodge numbers are
$$h^{1,1}(\Yf)=\left\{\begin{array}{rr}73&\mbox{ if }g=2\\56&\mbox{ if }g=3\end{array}\right.,\ \ \
h^{2,1}(\Yf)=\left\{\begin{array}{rr}1&\mbox{ if }g=2\\2&\mbox{ if
}g=3\end{array}\right..
$$
\end{prop}
\bprf
The existence of the desingularization of $\left(E_i\times \Sf\right)/\left(\alpha_E^3\times \alpha_S\right)$ is proved in \cite{cynk hulek} (indeed the action of $\alpha_E^3\times \alpha_S$ linearizes near the fixed locus as prescribed by the hypothesis of \cite[Proposition 4.1]{cynk hulek}). Here we describe explicitly the desingularization of this quotient in order to compute the Hodge numbers of $\Yf$.\\
{\it Step 1: the fixed locus.}
The fixed locus of $\alpha_S$ on the K3 surface $S_{f_g}$ consists of $c_g$ rational curves $C_i$, $i=1,\ldots, c_g$ and $2d_g$ points $R_{h,j}$, $h=1,2$, $j=1,\ldots d_g$. The fixed locus of $\alpha_S^2$ consists of course of the rational curves $C_i$ and of $d_g$ rational curves $D_j$, $j=1,\ldots d_g$. The points fixed for $\alpha_S$ are on the curves fixed for $\alpha_S^2$, in particular $R_{h,j}\in D_j$, $h=1,2$, $j=1,\ldots d_g$.\\
By the Proposition \ref{prop: non sympl autom on S}, we have $c_g=\left\{\begin{array}{ll} 3&\mbox{ if }g=2\\2&\mbox{ if }g=3\end{array}\right.,$  and $d_g=\left\{\begin{array}{ll} 5&\mbox{ if }g=2\\4&\mbox{ if }g=3\end{array}\right..$
The fixed locus of $\alpha_E^3$ on $E$ is made up of two 2-torsion points $P_1$, $P_2$ and the fixed locus for $\alpha_E^2$ consists of the 2-torsion group of $E$, i.e.\ of $E[2]=\{P_1, P_2, Q_1, Q_2\}$. We observe that $\alpha_E^3(Q_1)=Q_2$.\\
Hence the fixed locus of $\left(\alpha_E^3\times \alpha_S\right)^2$ consists of: \\ $\bullet$ $2c_g$ curves of type $\pc$ $i=1,2$, $j=1,\ldots c_g$,\\ $\bullet$ $2c_g$ curves of type $\qc$, $i=1,2$, $j=1,\ldots c_g$,\\ $\bullet$ $2d_g$ curves of type $\pd$ $i=1,2$, $j=1,\ldots d_g$,\\ $\bullet$ $2d_g$ curves of type $\qd$, $i=1,2$, $j=1,\ldots d_g$.\\The curves $\pc$ are in the fixed locus of $\alpha_E^3\times \alpha_S$, the curves  $\pd$ are invariant for $\alpha_E^3\times\alpha_S$, but they are not fixed: on each of them there are two fixed points. The other curves are switched by $\alpha_E^3\times \alpha_S$, indeed $(\alpha_E^3\times \alpha_S)(Q_1\times D_j)=Q_2\times D_j$ and $(\alpha_E^3\times \alpha_S)(Q_1\times C_j)=Q_2\times C_j$.\\

{\it Step 2: the desingularization $\Zf$ of the quotient $\left(E_i\times \Sf\right)/\left(\alpha_E^3\times \alpha_S\right)^2$.}
The curves fixed by $(\alpha_E^3\times \alpha_S)^2$ are $4c_g+4d_g$.
Let $\widetilde{E_i\times \Sf}$ be the blow up of $E_i\times \Sf$
along the fixed curves. The automorphism
$\alpha_E^3\times\alpha_S$ induces an automorphism
$\widetilde{\alpha_E^3\times \alpha_S}$ on $\widetilde{E_i\times
\Sf}$ of order 4, whose square fixes the exceptional divisors on the fixed curves.  The exceptional divisors on $\widetilde{E_i\times
\Sf}$ are $\mathbb{P}(\mathcal{N}_{(\pc)/(E_i\times \Sf)})$ (resp.\
$\mathbb{P}(\mathcal{N}_{(\pd)/(E_i\times \Sf)})$, $\mathbb{P}(\mathcal{N}_{(\qd)/(E_\times \Sf)})$,  $\mathbb{P}(\mathcal{N}_{(\qd)/(E_\times \Sf)})$), so they are the projectivization of a rank 2 bundle over $\pc$ (resp.\ $\pd$, $\qc$, $\qd$) which are rational curves. Hence the exceptional divisors are $\mathbb{P}^1$-bundles over rational curves.  We will denote the exceptional divisor on a curve $F$ as $\widetilde{F}$.\\ 
Let $\Zf$ be the quotient of $\widetilde{E_i\times \Sf}$ by $\widetilde{(\alpha_E^3\times\alpha_S)^2}$. This is a smooth 3-fold, and the quotient map $\widetilde{E_i\times \Sf}\ra \Zf$ is a 2:1 cover ramified on $\widetilde{\pc}$, $\widetilde{\pd}$, $\widetilde{\qc}$, $\widetilde{\qd}$. We have the following commutative diagram:
$$\begin{array}{cccr}E_i\times \Sf&\leftarrow&\widetilde{E_i\times \Sf}\\
\downarrow &&\downarrow\\
(E_i\times
\Sf)/\alpha_E^2\times\alpha_S^2&\leftarrow&\widetilde{(E_i\times
\Sf)/\alpha_E^2\times\alpha_S^2}\simeq \Zf&\end{array}$$
where $\widetilde{(E_i\times \Sf)/\alpha_E^2\times\alpha_S^2}$ is the desingularization the quotient of $(E_i\times \Sf)/\alpha_E^2\times\alpha_S^2$.\\

{\it Step 3: the desingularization of the quotient $\left(E_i\times \Sf\right)/\left(\alpha_E^3\times \alpha_S\right)$.} 
The automorphism $\alpha_E^3\times \alpha_S$ induces an automorphism of order 2, $\alpha_Z$, on $\Zf$. Since $(\alpha_E^3\times \alpha_S)(Q_1\times C_j)=Q_2\times C_j$ and  $(\alpha_E^3\times \alpha_S)(Q_1\times D_j)=Q_2\times D_j$, the divisors $\widetilde{\qc}$ and $\widetilde{\qd}$ are not invariant and in particular there are no elements in the fixed locus on these divisors. The divisors $\widetilde{\pc}$ and $\widetilde{\pd}$ are $\alpha_Z$-invariant. Since $\alpha_E^3\times \alpha_S$ fixes two points on $\pd$, $\alpha_Z$ fixes two fibers in the $\mathbb{P}^1$-bundle $\widetilde{\pd}$ (the fibers over the fixed points). These fibers are clearly two rational curves.\\
The involution $\alpha_Z$ fixes the bases of the $\mathbb{P}^1$-bundles $\widetilde{\pc}$ and acts non trivially on the fibers. Since the fibers are rational curves, $\alpha_Z$ fixes two points on each fiber, so it fixes two sections on $\widetilde{\pc}$. The base of the bundles $\widetilde{\pc}$ are the rational curves $\pc$, thus the sections of these bundles are rational curves. Hence $\alpha_Z$ fixes $4d_g+4c_g$ rational curves on $\Zf$. \\
Let $\widetilde{\Zf}$ be the blow up of $Z_{f_g}$ in the fixed locus of $\alpha_Z$, let $\widetilde{\alpha_Z}$ be the automorphism induced by $\alpha_Z$ on $\widetilde{Z_{f_g}}$. The fixed locus of $\widetilde{\alpha_Z}$ on $\widetilde{\Zf}$ is smooth of codimension 1, hence $\Yf:=\widetilde{\Zf}/\widetilde{\alpha_Z}$ is smooth. As before we have the following commutative diagram:
$$\begin{array}{cccr}\Zf&\leftarrow&\widetilde{\Zf}\\
\downarrow &&\downarrow\\
\Zf/\alpha_{\Zf}&\leftarrow&\widetilde{\Zf/\alpha_{Z}}\simeq \Yf&\end{array}$$ 
where $\widetilde{\Zf/\alpha_{Z}}$ is the desingularization the quotient of $\Zf/\alpha_Z$. Putting together the two quotients of order 2 one obtains:
$$\begin{array}{ccccccc}E_i\times \Sf&\leftarrow&\widetilde{E_i\times \Sf}\\
\downarrow &&\downarrow\\
(E_i\times \Sf)/\alpha_E^2\times\alpha_S^2&\leftarrow&\Zf&\leftarrow&\widetilde{\Zf}\\
&&\downarrow &&\downarrow\\
&&\Zf/\alpha_{\Zf}&\leftarrow&\widetilde{\Zf/\alpha_{Z}}\simeq
\Yf&\end{array}$$ which shows that $\Yf$ is a desingularization of
the quotient $E_i\times \Sf$ by $\alpha_E^3\times \alpha_S$.\\

{\it Step 4: the Hodge numbers $h^{2,1}$ and $h^{1,1}$.} Let $A$ be a K\"ahler 3-fold admitting a finite automorphism $\phi$ fixing a finite set of points and of disjoint curves. Let $\widetilde{A}$ be the blow up of $A$ along the fixed locus of $\phi$ and $\widetilde{\phi}$ be the automorphism induced on $\widetilde{A}$ by $\phi$. We can compute the Hodge structure of $\widetilde{A}$ by \cite[Th\'eor\`em 7.31]{voisin}. The $\widetilde{\phi}$-invariant part for the action of  $\widetilde{\phi}$ on the Hodge structure of $\widetilde{A}$ depends only on the $\phi$-invariant part of the Hodge structure of $A$ and on the Hodge structure of the fixed locus  (cf.\ \cite[Proposition 10.3.2]{Rohde thesis}). In particular, if the fixed locus of $\phi$ is made up of $r$ disjoint rational curves, one obtains:
\begin{equation}\label{fromula: hodge structure quotient}H^2(\widetilde{A},\Z)_1\simeq H^2(A,\Z)_1\oplus H^0(B,\Z)\
\ \ \mbox{ and }\ \ H^3(\widetilde{A},\Z)_1\simeq
H^3(A,\Z)_1\oplus H^1(B,\Z).\end{equation} 
Since $H^0(B, \Z)=\Z^r$ and $H^1(B,\Z)=0$ we obtain that $h^2(\widetilde{A}/\widetilde{\phi})=\dim(H^2(\widetilde{A},\Z)_1)=h^2(A)+r$ and $h^3(\widetilde{A}/\widetilde{\phi})=\dim(H^3(\widetilde{A},\Z)_1)=h^3(A).$
Recalling that the Hodge numbers $h^{i,0}$ are birational invariants one finds that $h^{1,1}(\widetilde{A}/\widetilde{\phi})=h^{1,1}(A)+r$ and $h^{2,1}(\widetilde{A}/\widetilde{\phi})=h^{2,1}(A)$.\\
We apply these results to $A=E_i\times \Sf$ and $\phi=(\alpha_E^3\times\alpha_S)^2$ and after that to $A=\Zf$ and $\phi=\alpha_Z$. So we obtain: 
$$h^{2,1}(\Zf)=dim\left(H^{2,1}(E_i\times \Sf)^{(\alpha_E^3\times\alpha_S)^2}\right)=g+1$$
because $H^{2,1}(E_i\times \Sf)^{(\alpha_E^3\times\alpha_S)^2}=H^{1,0}(E)\otimes H^{1,1}(S)_{i}\oplus H^{1,0}(E)\otimes H^{1,1}(S)_{-i}$ since the anti-invariant sublattice of $H^{1,1}(S)$ with respect to the action of $\alpha_S^2$ coincides with the direct sum of the eigenspaces for the eigenvalues $i$ and $-i$ for the action of $\alpha_S$;
$$h^{2,1}(\Yf)=dim\left(H^{2,1}(E_i\times \Sf)^{\alpha_E^3\times\alpha_S}\right)=g-1.$$
To construct $\Zf$ (the desingularization of the first quotient of order two) we blow up $4(c_g+d_g)$ rational curves, so 
$$h^{1,1}(\Zf)=dim\left(H^{1,1}(E_i\times \Sf)^{(\alpha_E^3\times\alpha_S)^2}\right)+4(c_g+d_g)=\left\{\begin{array}{ll} 51&\mbox{ if }g=2\\41&\mbox{ if }g=3\end{array}\right..$$ To construct $\Yf$ from the 3-fold $\Zf$ we blow up $4(c_g+d_g)$ rational curves, so we have $h^{1,1}(\Yf)=dim \left(H^{1,1}(\Zf)^{\alpha_Z}\right)+4(c_g+d_g)$.  Since $H^{1,1}(\Zf)^{\alpha_Z}$ decomposes in the sum of $H^{1,1}(E_i\times\Sf)^{\alpha_E^3\times\alpha_S}$ and of the space generated by $\widetilde{Q_1\times C_j}+\widetilde{Q_2\times C_j}$, $\widetilde{Q_1\times D_j}+\widetilde{Q_2\times D_j}$, $\widetilde{\pc}$ and  $\widetilde{\pd}$ one obtains:
\begin{equation}\label{eq: h11 y}h^{1,1}(\Yf)=dim(H^{1,1}(E_i\times S)^{\alpha_E^3\times \alpha_S})+3(c_g+d_g)+4(c_g+d_g)=\left\{\begin{array}{ll}73&\mbox{ if }g=2\\56&\mbox{ if }g=3\end{array}\right..\end{equation}\eprf

\begin{rem}{\rm The Hodge numbers of $\Zf$ were computed in \cite{voisin miroirs} and \cite{borcea}, indeed $\Zf$ is obtained by applying the Borcea--Voisin construction to the K3 surface $\Sf$ and the elliptic curve $E_i$. In particular $\Sf$ is associated to the invariant $(r,a,\delta)=(18,4,1)$ if $g=2$ and $(r,a,\delta)=(16,6,1)$ if $g=3$.}\end{rem}
\begin{rem}{\rm In the computation of the Hodge numbers (in particular of $h^{2,1}$) we used the fact that the fixed loci of $\alpha_S$ and of $\alpha_S^2$ do not contain curves of positive genus.}\end{rem}
\begin{rem}\label{rem: YSI}{\rm There exists a desingularization $Y_{SI}$ of $\left(E_i\times X_{SI}\right)/\left(\alpha_E^3\times \alpha_X\right)$ which is a Calabi--Yau. Its Hodge numbers are
$h^{1,1}(Y_{SI})=90$, $h^{2,1}(Y_{SI})=0$. This can be proved as in the Proposition \ref{prop: Hodge numbers CY}.}\end{rem}

\begin{prop}\label{prop: variation hodge structure}
The variation of the Hodge structures of weight 3 $\{H^3(\Yf,\Q)\}_{f_g}$ depends only on the variation of the Hodge structures $\{H^1(\Cf,\Q)\}_{f_g}$. More precisely, there exists an isomorphism $\phi:H^1(\Cf,\Q)\ra H^3(\Yf,\Q)$ such that $\phi(H^{1,0}(\Cf,\Q)_{-i})=H^{3,0}(\Yf,\Q)$ and $\phi(H^{0,1}(\Cf,\Q)_{-i})=H^{2,1}(\Yf,\Q)$.
\end{prop}
\bprf By the previous proposition, the Hodge structure of $H^3(\Yf,\Q)$ is the sub-Hodge structure of $H^3(E_i\times \Sf,\Q)$ which is invariant under the action of $\alpha_E^3\times \alpha_S$. By \eqref{formula: decomp. eigensp. H2(S)} and
\eqref{formula: H30 Y} we have
$$H^{3,0}(\Yf,\C)=\langle\eta_E\wedge\eta_S\rangle=\langle \eta_E\wedge \eta_E\wedge \omega_C\rangle= H^{1,0}(E_i) \otimes H^{1,0}(E_i)\otimes H^{1,0}(\Cf)_{-i}\simeq H^{1,0}(\Cf)_{-i},$$ 
since $E_i$ is rigid. Analogously by \eqref{formula: decomp. eigensp. H2(S)} and \eqref{formula: H21 Y} we have
$$H^{2,1}(\Yf)=H^{1,0}(E_i)\otimes \left(H^{1,1}(\Sf)\right)_i=H^{1,0}(E_i)\otimes H^{1,0}(E_i)\otimes H^{0,1}(\Cf)_i\simeq H^{0,1}(\Cf)_{-i}.$$\eprf
\begin{cor}\label{cor: moduli space} The complex varieties $\Cf$, $\Sf$ and $\Yf$ are parametrized by the same Shimura variety (cf. Proposition \ref{prop: moduli space Cf}).\end{cor}

\subsection{The desingularization $W_{b^2}$ of $E_i\times Q_{b^2}/\alpha_E^3\times \alpha_Q$.}\label{subsection: Wb}
In section \ref{subsection: isotrivial} we constructed a 1-dimensional family of K3 surfaces, $Q_{b^2}$, admitting a non symplectic automorphism of order 4, $\alpha_Q$, such that the fixed loci of $\alpha_Q$ and $\alpha_Q^2$ are made up of points and rational curves (cf.\ Proposition \ref{prop: Qb2}). Considering the quotient $E_i\times Q_{b^2}/\alpha_E^3\times \alpha_Q$ one can prove the existence of a smooth Calabi--Yau 3-fold, $W_{b^2}$, which desingularizes it, as in the Proposition \ref{prop: Hodge numbers CY}. Also the computation of the Hodge numbers of $W_{b^2}$ is similar to the one in proof of the Proposition \ref{prop: Hodge numbers CY}, the only difference being the fixed locus of $\alpha_E^3\times \alpha_S$: on the surface $Q_{b^2}$ there are two rational curves which are fixed by $\alpha_Q^2$, but which are switched by $\alpha_Q$ (this never happens on the surface $\Sf$). This does not change the computation of $h^{2,1}$, but has an effect on $h^{1,1}$, indeed if $Q$ is a K3 surface with a non symplectic automorphism of order 4, $\alpha_Q$, and there are $a$ curves fixed by $\alpha_Q^2$ but switched by $\alpha_Q$, $c$ curves fixed by $\alpha_Q$ and $d$ curves fixed by $\alpha_Q^2$ and invariant for $\alpha_Q$, then formula \eqref{eq: h11 y} becomes:
$$h^{1,1}(\widetilde{E_i\times Q/(\alpha_E^3\times \alpha_Q)})=dim \left(H^{1,1}(E_i\times S)^{\alpha_E^3\times \alpha_Q}\right)+3(c+d)+2a+4(c+d).$$
In our cases $a=2$, $c=2$, $d=4$, $dim \left(H^{1,1}(E_i\times S)^{\alpha_E^3\times \alpha_Q}\right)=15$, thus $$h^{2,1}(W_{b^2})=1,\ \ \ h^{1,1}(W_{b^2})=15+18+4+24=61.$$
As in Proposition \ref{prop: variation hodge structure} and Corollary \ref{cor: moduli space}, one finds that the variation of the Hodge structures of this family depends only on variation of the Hodge structures of an elliptic curve.\\
This family of Calabi--Yau is well known (cf.\ \cite{borcea}, \cite[Example 1]{Rohde}, \cite{voisin miroirs}) indeed $W_{b^2}$ is the desingularization of the quotient $E_i\times E_i\times E/\langle 1\times \alpha_E^2\times \iota_E,\alpha_E^3\times\alpha_E\times 1\rangle$ (cf.\ Remark \ref{rem: desingularization}) and hence  of $X_{SI}\times E/(\alpha_X^2\times \iota)$. Thus it is obtained by the Borcea--Voisin construction applied to the K3 surface $X_{SI}$ (which is the unique K3 surface admitting a non symplectic involution, $\alpha_X^2$, acting trivially on the N\'eron--Severi group and fixing 10 curves) and an elliptic curve.
\subsection{The desingularization $M$ of the quotient $E_i\times R/\alpha_E^3\times \alpha_R$}\label{subsection: 3 dim fam CY} Let the family $R$ of K3 surfaces be as in the Section \ref{prop: a 3 dimensional family}. 
Considering the quotient $E_i\times R/\alpha_E^3\times \alpha_R$ one can prove the existence of a smooth Calabi--Yau 3-fold, $M$, which desingularizes it, and compute the Hodge numbers of $M$ as in the Proposition \ref{prop: Hodge numbers CY}. In particular
$$h^{2,1}(M)=3,\ \ h^{1,1}(M)=11+3(1+3)+4(1+3)=39,$$
and the variation of the Hodge structures of $M$ depends only on the variation of the Hodge structures of the curve $B$ such that $\widetilde{B\times E_i/G}\simeq R$ (cf.\ Section \ref{section: geometrical description}).
 
\section{Picard--Fuchs equations of the 1-dimensional families}\label{section: PF}
\subsection{The family $\Yf$.}\label{subsect: PF Y} By the Proposition \ref{prop: variation hodge structure} it is clear that the Picard--Fuchs equations of the varieties $\Sf$ and $\Yf$ are given by the Picard--Fuchs equations of $\Cf$. In particular if $g=2$ we have a one dimensional family of Calabi--Yau 3-folds $\Yf$ (and of K3 surfaces $\Sf$). This situation is very similar to the one described in \cite[Section 2]{noi}. In this case one can explicitly compute the Picard--Fuchs equation of $\Cf$ (and hence of $\Sf$ and $\Yf$). Indeed (by the Section \ref{sec: the curves Cf}) one can assume that $C_{f_2}$ is given by the equation $$w^4=t(t-1)^2(t-\lambda)^2.$$ So the curves $\Cfd$ have an equation of type 
$$
C_\lambda\,:\; z^N\,=\,r^A(r-1)^B(r-\lambda)^C
$$
for some $\lambda\in \C$. In particular, $N=4$ and one can assume that $A=1$ and $B=C=2$. The holomorphic one forms on this curve are
$\omega(0,0,0;1)$ and $\omega(0,1,1;3)$ where
$$
\omega(\alpha,\beta,\gamma;l)\;:=\;
\frac{r^\alpha(r-1)^\beta(r-\lambda)^\gamma}{z^l}\rmd r.
$$
To find the Picard--Fuchs equations, let:
$$
a\,:=\,-\alpha+(lA/N),\qquad
b\,:=\,-\beta+(lB/N),\qquad
c\,:=\,-\gamma+(lC/N).
$$
An explicit computation shows that
$$
\left(\lambda(1-\lambda)\frac{\partial^2}{\partial \lambda^2}+
(a+c-\lambda(a+b+2c))\frac{\partial}{\partial \lambda}-c(a+b+c-1)\right)
\omega(\alpha,\beta,\gamma;l)\;=\;c\rmd h
$$
where $h$ is the rational function
$$
h\,:=\,\frac{r^{\alpha+1}(r-1)^{\beta+1}(r-\lambda)^{\gamma-1}}{z^l}\,=\,
r^{1-a}(r-1)^{1-b}(r-\lambda)^{-1-c}.
$$
In particular the Picard--Fuchs equation of $\Yfd$ is not a differential equation of order 4, but of order 2. This implies that the Picard--Fuchs equation of the family $\Yf$ has no a solution with maximal unipotent monodromy.
\subsection{The family $W_{b^2}$} In the Section \ref{subsection: Wb} we constructed another 1-dimensional family of Calabi--Yau 3-folds, $W_{b^2}$, and we proved that the variation of its Hodge structures depends only on the variation of the Hodge structures of an elliptic curve. So the Picard-Fuchs equation of $W_{b^2}$ is the one of the elliptic curve $z^2=r(r-1)(r-\lambda)$ (which is well known, and could be computed specilizing the general results presented in the Section \ref{subsect: PF Y}):
$$
\left(\lambda(1-\lambda)\frac{\partial^2}{\partial \lambda^2}+
(1-2\lambda)\frac{\partial}{\partial \lambda}-\frac{1}{4}\right)\frac{
dr}{z}\;=\;
\frac{1}{2}\rmd (r^{1/2}(r-1)^{1/2}(r-\lambda)^{-3/2}).
$$
As before this shows directly that there is no solution with maximal unipotent monodromy.

\subsection{The families $Y_{f_3}$ and $M$} 
More in general, the Calabi--Yau 3-folds obtained in the Section \ref{section: Calabi Yau}, as desingularization of $E_i\times S/G$ for a certain K3 surface $S$ and a certain group $G$, admit an automorphism of order 4, induced by $\alpha_E\times 1_S$, which acts as a multiplication by $i$ on the $H^{3,0}$ of the Calabi--Yau. Since this automorphism commutes with the differentiation on the parameters of the surface $S$, the Picard--Fuchs equations of these families of Calabi--Yau can not admit solutions with maximal unipotent monodromy (cf.\ \cite{Rohde}). Moreover the fact that $S$ is obtained as desingularization of $E_i\times B/G'$ for certain curve $B$  and group $G'$ implies that the Picard--Fuchs equations of the K3 surface and of the Calabi--Yau 3-fold depend only on the Picard--Fuchs equations of the curve $B$. Thus the Picard--Fuchs equations of the family $Y_{f_3}$ and of the family $M$ do not admit a solution with maximal unipotent monodromy and are the Picard--Fuchs equations of a curve.

\end{document}